\documentclass[10pt, leqno]{amsart}
\usepackage{amscd}
\usepackage{amssymb}
\usepackage{amsmath}
\usepackage{amsfonts}
\usepackage{stmaryrd}
\usepackage{amsthm}
\usepackage[cmtip,all]{xypic}
\usepackage{pifont}
\usepackage{faktor}

\usepackage[usenames, dvipsnames]{color}

\theoremstyle{plain}
\newtheorem{Lem}{Lemma}[section]

\newtheorem{Thm}[Lem]{Theorem}

{\theoremstyle{definition} 

\newtheorem{Rk}[Lem]{Remark}
\newtheorem{Def}[Lem]{Definition}}

\newcommand{\zig}{\addtocounter{Lem}{1}\tag{\theLem}}
\newcommand{\G}{\mathbb{Z}_p^\times}

\newcommand{\padics}{\mathbb{Z}_p}

\def\:{\colon}
\DeclareMathOperator*{\colim}{colim}
\DeclareMathOperator*{\holim}{holim}

\DeclareMathOperator*{\lims}{lim^\mathit{s}}
\DeclareMathOperator*{\limone}{lim^1}

\begin{document}
\title{A continuous $p$-adic action on the $K(2)$-local 
algebraic $K$-theory of $p$-adic complex $K$-theory}
\author{Daniel G. Davis}
\begin{abstract} 
Let $p$ be a prime, let $KU_p$ be $p$-complete complex $K$-theory, and let $\mathbb{Z}_p^\times$ denote the group of units in the $p$-adic integers. The $p$-adic Adams operations induce an action of the profinite group $\G$ on $KU_p$, and hence, on the algebraic $K$-theory spectrum $K(KU_p)$. For $p \geq 5$, we give an elementary construction of the continuous homotopy fixed point spectrum 
$(L_{K(2)}K(KU_p))^{hG}$, where $K(2)$ is the second Morava $K$-theory and $G$ is any closed subgroup of $\G$. Also, for each $G$, we show that there is an associated strongly convergent 
homotopy fixed point spectral sequence whose $E_2$-term is given by Jannsen's continuous group cohomology, with $E_2^{s,\ast} = 0$, for all $s > 2$. This work is related to a conjecture of Ausoni and Rognes. 
\end{abstract}

\maketitle
\section{Introduction}\label{sectionone}

Let $p$ be a prime, $\mathbb{Z}_p$ the $p$-adic integers, 
$L_{K(1)}S^0$ the Bousfield localization of the sphere spectrum with respect to 
the first Morava $K$-theory spectrum, and $KU_p$ $p$-complete complex $K$-theory. Thus, 
\[\pi_\ast(KU_p) = \mathbb{Z}_p[u^{\pm 1}],\] where $\pi_0(KU_p) = \mathbb{Z}_p$ 
and $|u| = 2$. Also, let $\mathbb{Z}_p^\times$ denote the group of units in 
$\mathbb{Z}_p$. By \cite{Pgg/Hop0, AndreQuillen}, there is a $\mathbb{Z}_p^\times$-action 
on the commutative $S^0$-algebra 
$KU_p$ via the $p$-adic Adams operations and this action is by maps of 
commutative $S^0$-algebras. Given a commutative $S^0$-algebra $A$, 
$K(A)$ denotes the algebraic $K$-theory spectrum of $A$, and $K(-)$ is a functor on the category of 
commutative $S^0$-algebras. It follows that $\mathbb{Z}_p^\times$ acts on 
$K(KU_p)$ by maps of commutative $S^0$-algebras.

For the rest of this paper, let $p \geq 5$ and let $V(1)$ be the type $2$ Smith-Toda 
complex $S^0/(p, v_1)$, a homotopy commutative ring spectrum. There is a $v_2$-self-map 
$v \: \Sigma^{2p^2-2} V(1) \to V(1)$ and $v_{2}^{-1}V(1)$ denotes the 
mapping telescope associated to $v$. 
In \cite[paragraph containing (0.1)]{acta}, \cite[Conjecture 4.2]{rognesguido}, 
and \cite[page 46; Remark 10.8]{jems},
Christian Ausoni and John Rognes conjectured that the unit map 
\[L_{K(1)}S^0 \to KU_p\] in the $K(1)$-local category induces a 
weak equivalence 
\begin{equation}\label{themap}\zig
K(L_{K(1)}S^0) \wedge v_2^{-1}V(1) \to K(KU_p)^{h\mathbb{Z}_p^\times} 
\wedge v_2^{-1}V(1),\end{equation} where the conjectural object $K(KU_p)^{h\mathbb{Z}_p^\times}$ is a 
continuous homotopy fixed point spectrum that is formed with respect to the conjecturally 
continuous action of the profinite group $\mathbb{Z}_p^\times$ on 
$K(KU_p)$. 

We point out that as \cite{acta, rognesguido, jems} explain, the above conjecture consists of 
$n = 1$ cases of a more general conjecture made by Ausoni and Rognes for every positive integer and every prime. There is no published construction of $K(KU_p)^{h\mathbb{Z}_p^\times}$ and, according to \cite[Remark 1.5]{padicspectra2}, the 
only models for it currently are a ``candidate definition" using condensed 
spectra in the sense of Clausen-Scholze and, possibly, a 
pyknotic version of this construction (in the framework of \cite{pyknotic}). As far as we know, documents about these models are not available. 

After some recollections in Section \ref{recall}, we give in Section \ref{construction} an elementary construction 
of 
$L_{K(2)}K(KU_p)$ as the homotopy limit of a tower of discrete $\G$-spectra. Here, $L_{K(2)}(-)$ denotes Bousfield localization with respect to the second Morava $K$-theory spectrum and we are using the term ``discrete $\G$-spectra" in the sense of 
\cite[Section 2.3]{joint}. In Definition \ref{keydef}, by following a familiar pattern in algebraic $K$-theory, we see that the aforementioned homotopy limit 
yields the continuous homotopy 
fixed point spectrum 
\[(L_{K(2)}K(KU_p))^{hG},\] where $G$ is any closed subgroup of 
\[\G \cong \padics \times C_{p-1}\] and $C_{p-1}$ is a cyclic group of order $p-1$. 

We also have a result about a homotopy fixed point spectral sequence for \[\pi_\ast\bigl((L_{K(2)}K(KU_p))^{hG}\bigr).\] 
To describe this, we need to make some preliminary remarks. 

We work in the stable model 
category $Sp^\Sigma$ of symmetric spectra of 
simplicial sets and we let 
\[(-)_f \: Sp^\Sigma \rightarrow 
Sp^\Sigma, \ \ \ Z \mapsto Z_f\] be functorial fibrant replacement: given a spectrum $Z$, 
there is a natural map $Z \rightarrow Z_f$ that is a trivial cofibration, 
with $Z_f$ 
fibrant. By \cite[Propositions 3.7, 5.1]{Mahowald/Sadofsky} and 
\cite[Section 2]{HoveyCech}, there is a tower
\begin{equation}\label{tower}\zig
V(1) =: M_0 \leftarrow M_1 \leftarrow \cdots \leftarrow M_i \leftarrow \cdots\end{equation} of finite type $2$ complexes having the following properties:
\begin{itemize}
\item
for any spectrum $Z$, there is an equivalence  
\[L_{V(1)}Z \simeq \holim_{i} (Z \wedge M_i)_f\,;\]
\item
each $M_i$ has 
enough algebraic structure in the stable homotopy category (see \cite{Devinatzrings} and \cite[Section 4]{HS} for the details) that $\pi_0(M_i)$ is a commutative ring and for every integer $t$ and any spectrum $Z$, $\pi_t(Z \wedge M_i)$ is a unitary 
module over $\pi_0(M_i)$; 
\item
there is an increasing sequence $m_0 = 1, m_1, m_2, ..., m_i, ...$ of integers, such that 
for each $i \geq 0$, there is an isomorphism $\pi_0(M_i) \cong \mathbb{Z}/{p^{m_i}\mathbb{Z}}$ of rings.
\end{itemize} 
In (\ref{tower}), after $M_0$, 
there are suspensions on each $M_i$ that -- as is common practice -- we have omitted. 

If $A$ is a discrete $G$-module, $H^\ast_c(G, A)$ denotes the continuous cohomology of $G$ with coefficients in $A$. Given a tower $\{A_i\}_{i \geq 0}$ in the category of discrete $G$-modules, 
$H^\ast_\mathrm{cont}(G; \{A_i\}_{i \geq 0})$ denotes continuous cohomology in the sense of Jannsen \cite{Jannsen}. We let $L_2(-)$ be Bousfield localization with respect to the coproduct $K(0) \vee K(1) \vee K(2)$ of Morava $K$-theories, where $K(0)$ is the Eilenberg-Mac Lane spectrum $H\mathbb{Q}$. Also, for each $i$, the $G$-action on $\pi_t(K(KU_p) \wedge L_2M_i)$ is the action induced by 
letting $G$ act -- on the level of spectra -- on $K(KU_p)$ only. 

\begin{Thm}\label{mainhfpsSS}
Let $p \geq 5$ and let $G$ be a closed subgroup of $\G$. There is a 
strongly convergent homotopy fixed point spectral sequence
\[E_2^{s,t} \cong H^s_\mathrm{cont}(G; \{\pi_t(K(KU_p) \wedge L_2M_i)\}_{i \geq 0}) \Longrightarrow \pi_{t-s}\bigl((L_{K(2)}K(KU_p))^{hG}\bigr),\] where for 
each $t \in \mathbb{Z}$,
\[E_2^{s,t} \cong \begin{cases}
\displaystyle{\bigl(\lim_i \pi_t(K(KU_p) \wedge L_2M_i)\bigr)^{\mspace{-3mu}G}},
&s = 0;\\
\displaystyle{\limone_i H^1_c(G, \pi_t(K(KU_p) \wedge L_2M_i))},
& s = 2;\\
0
& s \geq 3,\end{cases}\] and when $s = 1$ and $t \in \mathbb{Z}$, there is the short exact sequence
\[0 \to \limone_i \bigl(\pi_t(K(KU_p) \wedge L_2M_i)\bigr)^{\mspace{-2.5mu}G} 
\to E_2^{1,t} \to \lim_i H^1_c(G, \pi_t(K(KU_p) \wedge L_2M_i)) \to 0.\] 
\end{Thm} 

The proof of Theorem \ref{mainhfpsSS} is in Section \ref{hfpSS}. We are working on a result related to when $E_2^{2,\ast}$ in Theorem \ref{mainhfpsSS} vanishes, and when $G = \G$, we are working on expressing the above $E_2$-page in terms of $\pi_\ast(K(L_p))$, where $L_p$ is the $p$-complete Adams summand (as in \cite{acta}). 

We briefly discuss the input that is needed for our work. In Section \ref{recall}, we 
recall from \cite[Theorem 1.13]{padicspectra2} that if a $\G$-spectrum $X$ has degreewise 
finite homotopy groups, then it can be realized as a discrete $\G$-spectrum. Key input for 
this result is the fact that by 
\cite[Theorem 2.9]{symondsetal}, the natural map
\[H^\ast_c(\padics, A) \to H^\ast(\padics, A)\] is a graded isomorphism whenever $A$ is a finite discrete $\padics$-module. In Section \ref{construction}, we use the aforementioned realization result to show that for each $i \geq 0$, $K(KU_p) \wedge M_i$ is equivalent to a discrete 
$\G$-spectrum. These equivalences are possible due to an application of the thick subcategory theorem \cite[Theorem 7]{nilpotencetwo} and the fact that 
$K(KU_p) \wedge V(1)$ has degreewise finite homotopy groups. This fact is a 
consequence of \cite{AusoniInventiones} and also work in \cite{localization} and \cite{BDR}, 
\cite{asterisque}, and \cite{RognesNotes} -- the way these works tie together to give the 
finiteness result is outlined in \cite[pages 6--7]{padicspectra2}. 

We would like to note that besides the input described 
above, our constructions of $L_{K(2)}K(KU_p)$ -- with its action by $\G$ continuous, the continuous homotopy fixed points $(L_{K(2)}K(KU_p))^{hG}$, and its homotopy fixed point spectral sequence are low-tech. We hope that the uncomplicated nature of this work has some appeal to the reader.  

In \cite[page 2]{LeeLevy}, David Jongwon Lee and Ishan Levy state that work of 
Shay Ben Moshe, Shachar Carmeli, Tomer Schlank, and Lior Yanovski ``combined with" \cite[Theorem 1.3.6]{HRW} by Jeremy Hahn, Arpon Raksit, and Dylan Wilson ``shows that $L_{K(2)}K(L_{K(1)}\mathbb{S}) \to L_{K(2)}K(KU_p)^{h\G}$ is an equivalence." Thus, our understanding is that the argument/work referred to by Lee and Levy involves a construction of 
$(L_{K(2)}K(KU_p))^{h\G}$. We do not know any of the details of this construction and so we do not know its relationship to the model for $(L_{K(2)}K(KU_p))^{h\G}$ given in the present paper. 

We make a few comments about our notation. Given $t \in \mathbb{Z}$ and a spectrum $Z$, $\pi_t(Z)$ denotes $[S^t, Z]$, the 
set of morphisms $S^t \to Z$ in the homotopy category $\mathrm{Ho}(Sp^\Sigma)$, 
where here, $S^t$ is a fixed cofibrant and fibrant model 
for the $t$-th suspension of the sphere spectrum. By ``$\holim$," we mean the homotopy limit for $Sp^\Sigma$, as defined in 
\cite[Definition 18.1.8]{hirschhorn}. 

\subsection*{Acknowledgements} I thank John Rognes for helpful tips about the degreewise finiteness of each $\pi_\ast(K(KU_p) \wedge M_i)$. 

\section{A discrete $H$-spectrum naturally associated to an $H$-spectrum, when $H$ is profinite}\label{recall}

Let $H$ be any profinite group and let $Z$ be any spectrum. 
We define $\mathrm{Sets}(H,Z)$ to be the $H$-spectrum whose $k$th pointed simplicial set 
$\mathrm{Sets}(H,Z)_k$ has $l$-simplices 
\[\mathrm{Sets}(H,Z)_{k,l} := \mathrm{Sets}(H,Z_{k,l}),\] the $H$-set of all functions 
$H \to Z_{k,l}$, for all $k, l \geq 0$. The $H$-action 
on $\mathrm{Sets}(H,Z_{k,l})$ is given by 
\[(h \cdot f)(h') = f(h'h), \ \ \ f \in \mathrm{Sets}(H,Z_{k,l}), \ h, h' \in H.\]
By \cite[Section 2]{padicspectra2}, when $X$ is an $H$-spectrum, there is a 
cosimplicial $H$-spectrum \[\mathrm{Sets}(H^{\bullet + 1}, X),\] where for 
$n \geq 0$, $\mathrm{Sets}(H^{\bullet + 1}, X)([n])$ is 
equal to the application of $\mathrm{Sets}(H,-)$ iteratively $n+1$ times to $X$. 

It is useful to note that if $K$ is an abstract group and $Y$ is a $K$-spectrum, then by functoriality, $Y_f$ is also 
a $K$-spectrum and the trivial cofibration 
$Y \rightarrow Y_f$ is $K$-equivariant. The following definition is based on 
\cite[Definition 4.4]{padicspectra2}, but is slightly more general. 

\begin{Def}\label{xdisO}
Let $X$ be an $H$-spectrum and let $\mathcal{O} = \{N_\lambda\}_{\lambda \in \Lambda}$ be an inverse 
system of open normal subgroups of $H$ ordered by inclusion, 
over $\Lambda$, a directed poset. We let
\begin{align*} 
X^\mathrm{dis}_\mathcal{O} := \colim_{\lambda \in \Lambda} 
\holim_\Delta \mathrm{Sets}(H^{\bullet + 1}, X_f)^{N_\lambda},\end{align*}
where the colimit is formed 
in $Sp^\Sigma$ and $(-)^{N_\lambda}$ denotes $N_\lambda$-fixed points. For each $\lambda$, the finite discrete group $H/N_\lambda$ acts on $\displaystyle{\holim_\Delta \mathrm{Sets}(H^{\bullet + 1}, X_f)^{N_\lambda}}$, so that the projection $H \to H/N_\lambda$ makes it a discrete 
$H$-spectrum. By \cite[Section 3.4]{joint}, colimits in the category of discrete $H$-spectra are formed in spectra, so we see that 
$X^\mathrm{dis}_\mathcal{O}$ is a discrete $H$-spectrum.  
\end{Def} 

By \cite[Lemma 4.7, proof of Theorem 4.9]{padicspectra2}, 
for any $H$-spectrum $X$, there is a zigzag 
\[\underbrace{X \xrightarrow{\,\simeq\,} X_f \xrightarrow{\,\simeq\,} \holim_\Delta \mathrm{Sets}(H^{\bullet+1}, X_f)}_{\text{the $H$-equivariant map} \ i_X} \xleftarrow{\phi_X} X^\mathrm{dis}_\mathcal{O}\] of $H$-equivariant maps, where $i_X$ is a weak equivalence of spectra and $\phi_X$ 
is induced by the inclusions $\mathrm{Sets}(H^{\bullet+1}, X_f)^{N_\lambda} \to 
\mathrm{Sets}(H^{\bullet+1}, X_f)$ of cosimplicial spectra.  

Now let $H = \G$ and let
\[\mathcal{O} = \{p^j\mathbb{Z}_p\}_{j \geq 0},\] where each $p^j\mathbb{Z}_p$ is the open normal subgroup of $\G$ that corresponds to the open normal subgroup $(p^j\mathbb{Z}_p) \times \{e\}$ of $\padics \times C_{p-1}$. By \cite[Theorem 1.13]{padicspectra2}, 
if the $\G$-spectrum 
$X$ has degreewise finite homotopy groups (that is, for each $t \in \mathbb{Z}$, the abelian group $\pi_t(X)$ is finite), then the 
map $\phi_X$ above is a weak equivalence in $Sp^\Sigma$. In this case, it is natural to identify the $\G$-spectrum $X$ 
with the discrete $\G$-spectrum $X^\mathrm{dis}_\mathcal{O}$, and hence, for any closed subgroup $G$ of $\G$, there is the 
continuous homotopy fixed point spectrum 
\[X^{hG} := (X^\mathrm{dis}_\mathcal{O})^{hG}.\] 

We recall the general definition of the above continuous homotopy fixed point spectrum, where 
$H$ is now an arbitrary profinite group. Let $(-)_{fH}$ denote functorial fibrant replacement for the model category of 
discrete $H$-spectra \cite[Theorem 2.3.2]{joint}. Unpacking the meaning of this yields that if 
$Y$ is any discrete $H$-spectrum, then there is a natural map 
$Y \to Y_{fH}$ that is $H$-equivariant, a trivial cofibration of spectra, and has a target that is a fibrant 
discrete $H$-spectrum. Then as in \cite[Section 3.1]{joint}, $Y^{hH}$ is obtained by taking the right derived functor of $H$-fixed points: 
\[Y^{hH} := (Y_{fH})^H.\] 


\section{The construction of $(L_{K(2)}K(KU_p))^{hG}$}\label{construction}

As recalled in Section \ref{sectionone}, $K(KU_p) \wedge V(1)$ has degreewise finite homotopy groups. It is straightforward to see that the full subcategory of the homotopy category of $p$-local finite spectra that consists of spectra $Z$ with the property that $K(KU_p) \wedge Z$ has degreewise finite homotopy groups is a thick subcategory. By the thick subcategory theorem \cite[Theorem 7]{nilpotencetwo}, since $V(1)$ is in this thick subcategory, it follows that any $K(1)$-acyclic $p$-local finite spectrum is in this subcategory, and hence, $K(KU_p) \wedge M_i$ is too, for all $i \geq 0$. 

For each $i \geq 0$, we regard $K(KU_p) \wedge M_i$ as a $\G$-spectrum by letting $\G$ act only on $K(KU_p)$. Since each $K(KU_p) \wedge M_i$ has degreewise finite homotopy groups, it follows from Section \ref{recall} that there is a zigzag
\[\bigl\{K(KU_p) \wedge M_i\bigr\}_{\mspace{-3mu}i} \mspace{-2mu}\xrightarrow{\,\simeq\,}\mspace{-2mu}
\bigl\{\holim_\Delta \mathrm{Sets}(H^{\bullet+1}, \mspace{-2mu}(K(KU_p) \wedge M_i)_f)\bigr\}_{\mspace{-3mu}i}
\mspace{-2mu}\xleftarrow{\,\simeq\,} \mspace{-2mu}\bigl\{(K(KU_p) \wedge M_i)^\mathrm{dis}_\mathcal{O}\bigr\}_{\mspace{-3mu}i}\] of morphisms of towers of $\G$-spectra, with each morphism a 
levelwise weak equivalence of spectra. Also, the rightmost tower is a tower in the category of 
discrete $\G$-spectra.

Let $H$ be any profinite group. We write ``$N \vartriangleleft_o H$" to signify that $N$ is an open normal subgroup of $H$. Notice that if $Y$ is a discrete $H$-spectrum and $Z$ is 
any spectrum with trivial $H$-action, then the induced $H$-action on 
\[Y \wedge Z \cong (\colim_{N \vartriangleleft_o H} Y^N) \wedge Z \cong 
\colim_{N \vartriangleleft_o H} (Y^N \wedge Z)\] makes $Y \wedge Z$ a discrete $H$-spectrum, 
since the projection $H \to H/N$ makes the $H/N$-spectrum $Y^N \wedge Z$ into a discrete $H$-spectrum. 

Now we make use of the 
above zigzag of towers. It is useful to note that $L_2(-)$ is a smashing localization: for any spectrum $Z$, $L_2Z \simeq Z \wedge L_2S^0$. Then by \cite[Remark 3.6]{DevFields} and \cite[Proposition 7.10]{HS}, there is an equivalence
\[L_{K(2)}K(KU_p) \simeq \holim_i (K(KU_p) \wedge M_i \wedge L_2S^0)_f.\]
By \cite[Corollary 5.3.3]{joint}, a fibrant discrete $\G$-spectrum is fibrant as a spectrum. 
Also, for each $i$, $(K(KU_p) \wedge M_i)^\mathrm{dis}_\mathcal{O} \wedge L_2S^0$ is a discrete 
$\G$-spectrum, where the $\G$-action is induced by the $\G$-action on $K(KU_p)$ and 
the trivial $\G$-action on each of $M_i$ and $L_2S^0$. 

We put the preceding facts together to 
conclude that there is an 
equivalence
\[L_{K(2)}K(KU_p) \simeq \holim_i ((K(KU_p) \wedge M_i)^\mathrm{dis}_\mathcal{O} \wedge L_2S^0)_{\mspace{-0mu}f\G}.\] In this equivalence, the right-hand side is the homotopy limit of a tower of 
discrete $\G$-spectra. Thus, it is natural to make the following definition.

\begin{Def}\label{keydef}
Let $G$ be a closed subgroup of $\G$. We 
define the continuous homotopy fixed point spectrum
\[(L_{K(2)}K(KU_p))^{hG} := \bigl(\holim_i ((K(KU_p) \wedge M_i)^\mathrm{dis}_\mathcal{O} \wedge L_2S^0)_{\mspace{-0mu}f\G}\bigr)^{hG},\] where 
\begin{align*}
\bigl(\holim_i ((K(KU_p) &\wedge M_i)^\mathrm{dis}_\mathcal{O} \wedge L_2S^0)_{\mspace{-1.5mu}f\G}\bigr)^{hG}\\
& := \mspace{-2mu} \holim_i \bigl(((K(KU_p) \wedge M_i)^\mathrm{dis}_\mathcal{O} \wedge L_2S^0)_{\mspace{-0mu}f\G}\bigr)^{hG}\mspace{-2mu}\end{align*} and for each $i$,
$\bigl(((K(KU_p) \wedge M_i)^\mathrm{dis}_\mathcal{O} \wedge L_2S^0)_{\mspace{-0mu}f\G}\bigr)^{hG}$ is 
the continuous $G$-homotopy fixed points of the given discrete $G$-spectrum (here, we are regarding the stated discrete $\G$-spectrum as a discrete $G$-spectrum). 
\end{Def}

\begin{Rk} The above definition follows a well-known script in algebraic $K$-theory: for example, 
see \cite[Proposition 3.1.2, last paragraph of 3.1, proof of Theorem 4.2.6]{Geisser}. For a general cofiltered diagram of discrete $G$-spectra, the script is expressed formally in \cite[Sections 4.4, 4.5]{joint}. In the case of towers -- as in Definition \ref{keydef}, \cite[Lemma 8.3, Remark 8.4]{cts} shows in the setting of Bousfield-Friedlander spectra that for a tower $\{X_i\}_i$ of discrete $G$-spectra, with each $X_i$ a fibrant spectrum, the continuous homotopy fixed point spectrum
\[(\holim_i X_i)^{hG} := \holim_i (X_i)^{hG}\] 
is the output of the right derived functor of 
the functor $\lim_i (-)^G$ from the category of towers of discrete $G$-spectra, equipped with the injective model structure, to spectra, when applied to the tower $\{X_i\}_i$.
\end{Rk}

\begin{Rk}
Let $H$ be any profinite group. There is the commutative diagram
\[\xymatrix@-.4pc{\{Y_i\}_i \ar[r] \ar[d] & \{Y'_i\}_i
\ar[d] \\ 
\{(Y_i)_{fH}\}_i \ar[r] & \{(Y'_i)_{fH}\}_i}\] in the category of towers of discrete $H$-spectra (in $Sp^\Sigma$, the default category of spectra in this paper), in which the top horizontal morphism is an arbitrary morphism in the category that levelwise consists of weak equivalences of spectra and the other three morphisms are given 
by fibrant replacement. Thus, the two vertical maps consist of levelwise weak equivalences of spectra. Since a morphism of discrete $H$-spectra that is a weak equivalence of spectra is a weak equivalence of discrete $H$-spectra, it follows that the bottom horizontal map consists of levelwise 
weak equivalences between fibrant objects in the model category of discrete $H$-spectra. Application of the right Quillen functor $(-)^H$ to this bottom map yields the morphism 
$\{(Y_i)^{hH}\}_i \xrightarrow{\,\simeq\,} \{(Y'_i)^{hH}\}_i$ of towers of spectra that consists of 
levelwise weak equivalences between fibrant spectra. Thus, there is a weak equivalence
\[\holim_i (Y_i)^{hH} \xrightarrow{\,\simeq\,}  \holim_i (Y'_i)^{hH}.\] This observation implies that for any $G$ closed in $\G$, there is the weak equivalence
\[\holim_i ((K(KU_p) \wedge M_i)^\mathrm{dis}_\mathcal{O} \wedge L_2S^0)^{hG} 
\xrightarrow{\,\simeq\,} \holim_i \bigl(((K(KU_p) \wedge M_i)^\mathrm{dis}_\mathcal{O} \wedge L_2S^0)_{f\G}\bigr)^{hG}.\] The target of this equivalence 
appears in Definition \ref{keydef}, and hence, the source of this equivalence -- which is a slightly 
simpler expression than the target -- can also 
be taken as the definition of 
$(L_{K(2)}K(KU_p))^{hG}$.   
\end{Rk}

%

\section{The homotopy fixed point spectral sequence for $(L_{K(2)}K(KU_p))^{hG}$}\label{hfpSS}

In this section, we prove Theorem \ref{mainhfpsSS}. We let $G$ be a closed subgroup of $\G$. 
If $H$ is a profinite group and $U$ is an open subgroup of $H$ such that $H^s_c(U, A) = 0$ for all discrete $U$-modules $A$, 
whenever $s > s_0$, for some integer $s_0$, then we say that $H$ has {\em finite virtual cohomological dimension}. Also, as is standard, $\mathrm{cd}_p(H)$ denotes the 
cohomological $p$-dimension of $H$.

The starting point is that $G$ has finite virtual cohomological dimension: the details for this assertion are, for example, in \cite[page 330]{cts} and the salient fact is that $\mathrm{cd}_p(\padics) = 1$. 
By \cite[Section 3.2, Theorem 3.2.1]{joint}, if $Y$ is a discrete $G$-spectrum that is fibrant as a 
spectrum, then there is a cosimplicial spectrum 
$\mathrm{Map}^c(G^\bullet, Y)$ that has the following properties: 
\begin{itemize}
\item
there is an equivalence $Y^{hG} \simeq \displaystyle{\holim_\Delta \mathrm{Map}^c(G^\bullet, Y)}$;
\item 
by \cite[Lemma 2.4.2, Section 4.6]{joint}, in cosimplicial degree $n$, where $n \geq 0$, there is an isomorphism
\[\pi_t\bigl(\mathrm{Map}^c(G^\bullet, Y)([n])\bigr) \cong \mathrm{Map}^c(G^n, \pi_t(Y)), \ \ \ t \in \mathbb{Z},\] where the right-hand side is the abelian group of continuous functions from the 
$n$-fold product $G^n$ to the discrete abelian group $\pi_t(Y)$ (if $n = 0$, then $G^n = \{e\}$, the trivial group);
\item
by \cite[Section 4.6]{joint}, there is a graded isomorphism \[H^\ast\bigl[\pi_t(\mathrm{Map}^c(G^\ast, Y))\bigr] \cong H^\ast_c(G, \pi_t(Y)), \ \ \ t \in \mathbb{Z},\] where the left-hand side is the cohomology of the cochain complex associated to the 
cosimplicial abelian group $\pi_t(\mathrm{Map}^c(G^\bullet, Y))$.
\end{itemize}

To simplify our notation, let
\[\mathbb{K}^\mathrm{dis}_i := ((K(KU_p) \wedge M_i)^\mathrm{dis}_\mathcal{O} \wedge L_2S^0)_{f\G}, \ \ \ i \geq 0.\] Also, we set 
\[\{i\} : = \{0 \shortrightarrow 1 \shortrightarrow 2 \shortrightarrow \cdots \shortrightarrow i \shortrightarrow \cdots\}^\mathrm{op},\] the opposite category of the 
directed poset $\{0 \shortrightarrow 1 \shortrightarrow \cdots \shortrightarrow i \shortrightarrow \cdots\}$. 
Then by applying the above recollections to Definition \ref{keydef}, we find that
\begin{equation}\label{doublyindexed}\zig
(L_{K(2)}K(KU_p))^{hG} \simeq \holim_i \holim_\Delta \mathrm{Map}^c(G^\bullet, \mathbb{K}^\mathrm{dis}_i) \cong \holim_{\Delta \times \{i\}} \mathrm{Map}^c(G^\bullet, \mathbb{K}^\mathrm{dis}_i)\end{equation} and there is the conditionally convergent homotopy 
spectral sequence 
\[E_2^{s,t}
\Longrightarrow \pi_{t-s}\Bigl(\,\holim_{\Delta \times \{i\}} \mathrm{Map}^c(G^\bullet, \mathbb{K}^\mathrm{dis}_i)\Bigr),\]
where 
\[E_2^{s,t} \cong \lims_{\Delta \times \{i\}} \pi_t\bigl(\mathrm{Map}^c(G^\bullet, \mathbb{K}^\mathrm{dis}_i)\bigr) \cong  \lims_{\Delta \times \{i\}} \mathrm{Map}^c(G^\bullet, \pi_t(\mathbb{K}^\mathrm{dis}_i))\] and the diagram $\mathrm{Map}^c(G^\bullet, \pi_t(\mathbb{K}^\mathrm{dis}_i))$ is 
the one obtained from isomorphisms stated above. 
By \cite[Proposition 3.1.2]{Geisser} (the particular case that we have here follows from 
\cite[Theorem 8.5 and the paragraph after it]{cts}, which is based on \cite{Geisser}, and 
\cite[proof of Theorem 3.2.1]{joint}) and (\ref{doublyindexed}), the above homotopy spectral 
sequence is the homotopy fixed point spectral sequence
\[E_2^{s,t} \cong H^s_\mathrm{cont}(G; \{\pi_t(\mathbb{K}^\mathrm{dis}_i)\}_{i \geq 0}) 
\Longrightarrow \pi_{t-s}\bigl((L_{K(2)}K(KU_p))^{hG}\bigr).\]
For each $i$, there are $\G$-equivariant isomorphisms
\[\pi_\ast(\mathbb{K}^\mathrm{dis}_i) \cong \pi_\ast(K(KU_p) \wedge M_i \wedge L_2S^0) 
\cong \pi_\ast(K(KU_p) \wedge L_2M_i),\] where in the last expression, $\G$ acts trivially on $L_2M_i$. These isomorphisms are natural in $i$, so that 
\[E_2^{s,t} 
\cong H^s_\mathrm{cont}(G; \{\pi_t(K(KU_p) \wedge L_2M_i)\}_{i \geq 0}), \ \ \ s \geq 0, \ t \in \mathbb{Z}.\]

For each $i \geq 0$, let 
\[\mathbb{K}_i := K(KU_p) \wedge L_2M_i.\]
By \cite[(2.1)]{Jannsen}, for each $s \geq 0$ and every $t \in \mathbb{Z}$, there is a short exact sequence 
\[0 \to \limone_i H^{s-1}_c(G, \pi_t(\mathbb{K}_i)) 
\to E_2^{s,t} \to \lim_i H^s_c(G, \pi_t(\mathbb{K}_i)) \to 0,\] where $H^{-1}_c(G, \pi_t(\mathbb{K}_i)) = 0$. Since $p \geq 5$, $\padics$ is the $p$-Sylow subgroup of $\G$. Thus,
\[\mathrm{cd}_p(G) \leq \mathrm{cd}_p(\G) = \mathrm{cd}_p(\padics) = 1,\] so that 
\[H^s_c(G, A) = 0, \ \ \ s \geq 2,\] whenever $A$ is a discrete $G$-module that as an abelian group is $p$-primary torsion. Since every \[\pi_t(\mathbb{K}_i) \cong \pi_t(K(KU_p) \wedge L_2S^0 \wedge M_i)\] 
is a unitary 
$(\mathbb{Z}/p^{m_i}\mathbb{Z})$-module,
\[H^s_c(G, \pi_t(\mathbb{K}_i)) = 0, \ \ \ s \geq 2, \ t \in \mathbb{Z}, \ i \geq 0.\] 
This last fact, together with the above short exact sequence, yields the remaining assertions in 
Theorem \ref{mainhfpsSS}.

\end{document}